\documentclass[epic,eepic,epsfig,envcountresetchap,envcountsame,draft]{amsart} 
\usepackage{amssymb}
\usepackage{amsmath}
\usepackage{amscd}
\usepackage{epic,eepic,epsfig} 
\makeindex  

\newcommand{\bC}{{\mathbb C}}

\newcommand{\bN}{{\mathbb N}}
\newcommand{\bP}{{\mathbb P}}
\newcommand{\bQ}{{\mathbb Q}}

\newcommand{\bZ}{{\mathbb Z}}
\newcommand{\Zb}{{\mathbb Z}}

\newcommand{\Gal}{{\rm Gal}}

\def\iii{{\Ic^\prime _2(X,Y)}}
 
\def\Cbb{{ C}}

\def\Nbb{{ N}}
\def\Pbb{{P}}
\def\Qbb{{ Q}}

\def\Zbb{{ Z}}

\def\Cc{{\mathcal C}}
\def\Gc{{\mathcal G}}
\def\Ic{{\mathcal I}}

\def\Lc{{\mathcal L}}
\def\Mc{{\mathcal M}}

\def\Uc{{\mathcal U}}
\def\Vc{{\mathcal V}}

\def\tfk{{\frak t}}

\def\01{{\overrightarrow{01}}}
\def\10{{\overrightarrow{10}}}
\def\2{{\overrightarrow{1\infty}}}

\def\deg{{\rm deg}}

\def\hpb{\hfill $\Box$}

\def\lra{\longrightarrow}

\def\al{\alpha}
\def\be{\beta}
\def\de{\delta}
\def\ga{\gamma}
\def\si{\sigma}
\def\la{\lambda}
\def\lam{\lambda}
\def\La{\Lambda}
\def\Lam{\Lambda}
\def\Ga{\Gamma}

\def\mod{{\rm mod}}

\def\log{{\rm log}}

\def\pd1{{\partial \Delta [1]}} 
\def\d1{{\Delta [1]}}

\def\zl{\Zbb _\ell}

\def\zlt{\Zbb _\ell^\times}

\def\Cbb{{\mathbb C}}

\def\Nbb{{\mathbb N}}
\def\Pbb{{\mathbb P}}
\def\Qbb{{\mathbb Q}}

\def\Ybb{{\mathbb Y}}
\def\Zbb{{\mathbb Z}}

\def\Cc{{\mathcal C}}

\def\Gc{{\mathcal G}}

\def\Ic{{\mathcal I}}

\def\Kc{{\mathcal K}}
\def\Lc{{\mathcal L}}
\def\Mc{{\mathcal M}}

\def\Uc{{\mathcal U}}
\def\Vc{{\mathcal V}}

\def\Yc{{\mathcal Y}}
\def\Zc{{\mathcal Z}}

\def\z2{{{\Zbb [{\frac{1}{2}}]}}}

\def\dfk{{\mathfrak d}}
\def\ffk{{\mathfrak f}}

\def\tfk{{\mathfrak t}}
\def\zfk{{\mathfrak z}}

\def\al{{\alpha}}
\def\be{{\beta}}
\def\ga{{\gamma}}
\def\si{{\sigma}}
\def\de{{\delta}}
\def\ka{{\kappa}}

\def\De{{\Delta}}

\def\01{{\overset{\to}{01}}}
\def\10{\overset{\to}{10}}
\def\2{\overset{\to}{1\infty}}

\def\hpb{\hfill $\Box$}

\def\lra{\longrightarrow}

\title{A weak Euler formula for l-adic Galois double zeta values}
\author{Zdzis{\l}aw Wojtkowiak}
\begin{document}

\date{\today} 

\begin{abstract}
The fact that the double zeta values $\zeta (n,m)$ can be written as a sum of products of two zeta values $\zeta (a )\zeta (a)$ and of $\zeta (n+m)$, 
whenever $n+m$ is odd is due to Euler. We shall show the weak version of this result for the l-adic Galois realization.
\end{abstract}

\maketitle 

\tableofcontents

\section{Introduction.}

The fact that the double zeta values 
\[
 \zeta (a+1, b+1)=\int _0^1 {\frac{dz}{1-z}},({\frac{dz}{ z}})^a,{\frac{dz}{1-z}},({\frac{dz}{z}})^b=\sum _{i_2>i_1=1}{\frac{1}{i_2^{b+1}}}
 {\frac{1}{i_1^{a+1}}}
\]
can be written as a sum of products of two zeta values $\zeta (\al )\zeta (\be )$ with $\al$ and $\be $ positive integers and of $\zeta (a+1+b+1)$, 
whenever $a+1+b+1$ is a positive odd integer is attributed to Euler.

\smallskip

The purpose of this paper is to prove a weak analogue of this result for l-adic Galois multi-zeta values.
Let us explain what we mean by this statement.
The Galois group
\[
 G_\bQ:= \Gal (\bar \bQ/\bQ)
\]
acts on the \'etale fundamental group
\[
 \pi _1(\bP ^1_{\bar \bQ}\setminus \{0,1,\infty \},\01 )\,,
\]
(see \cite{D} and \cite{I}).
Let us fix an embedding of $\bar \bQ $ into $\bC$. Let $p$ be a fixed rational prime number.
Let $x$ and $y$ be generators of $\pi _1(\bP ^1_{\bar \bQ}\setminus \{0,1,\infty \},\01 )$ as in \cite[page 154]{W1}.
Let
\[
 E:\pi _1(\bP ^1_{\bar \bQ}\setminus \{0,1,\infty \},\01 )\to \bQ_p\{\{X,Y\}\}
\]
be a continuous multiplicative map defined by 
\[
 E(x):=\exp (X)\;\;{\rm and}\;\;E(y):=\exp (Y)\,.
\]
Let $\pi$ be the canonical path from $\01$ to $\10$ on $\bP ^1_{\bar \bQ}\setminus \{0,1,\infty \}$, the interval $[0,1]$.

For $\si \in G_{\bQ}$, let us define a power series
\[
 \Lam _\pi (\si ):=E(\pi ^{-1}\cdot \si (\pi ))\in  \bQ_p\{\{X,Y\}\}\,.
\]
The coefficients of the power series $\Lam _\pi$ or $\log \Lam _\pi$, considered as functions on $G_\bQ$, are analogues of multi-zeta numbers 
(see \cite{W1}, where we pointed this analogy).
Let us denote by 
\[
 \la _{YX^aYX^b}
\]
the coefficient at $YX^aYX^b$ of the power series  $\La _\pi$. This coefficient we view as an analogue of the multi-zeta value $\zeta (a+1,b+1)$. 
The analogue of the Euler result mentioned at the beginning will be the following conjecture, 
which we state only for elements of $G_{\bQ (\mu _{p^\infty})}$.

\medskip

\noindent{\bf Conjecture A.}
If $a+1+b+1$ is odd then  $\la _{YX^aYX^b}$ is a linear combination of products $ \la _{YX^\al }\cdot \la _{YX^\be}$
with $\al +\be \leq a+b$ and of $ \la _{YX^\al }$ with $\al\leq a+b+1$.

\medskip

The actual result we shall prove is much weaker, hence ``weak Euler formula'' in the title of the paper. Let us define a subfield of $\bar \bQ$,
\[
 \Kc _1:=\bQ (\mu _{p^\infty})((1-\xi ^i_{p^n})^{\frac{1}{p^m}}\mid n,m \in \bN,\, 0< i<p^n)\,.
\]
We shall prove the following result as well as 
its generalization.
\medskip

\noindent{\bf Theorem B.} Let us assume that $a+1+b+1$ is odd. If $\si \in G_{\Kc _1}$ then $\la _{YX^aYX^b}(\si )=0$.

It seems clear that Conjecture A and Theorem B and its generalizations can be proved using the Drinfeld-Ihara-Deligne relations:
the $\bZ/2\bZ$-relation $\Lam _\pi (X,Y)\cdot \Lam _\pi (Y,X)=1$, the  $\bZ/3\bZ$-relation and the  $\bZ/5\bZ$-relation (see \cite{I1}), 
in the same way as in \cite{W7} we have calculated the coefficients $\lam _{YX^{2n-1}}$.

\medskip

The proof we present in this paper is however different. In \cite{W9} we have shown that the coefficient
\[
\la _{YX^{a_1} YX^{a_2}\ldots YX^{a_r} }
\]
of $\Lam _\pi $ at $YX^{a_1} YX^{a_2}\ldots YX^{a_r}$ can be express as an integral on $(\bZ _p)^r$ against the measure which we denoted by $G_r(\10 )$. 
We shall show that the measures
$G_r(\10 )$ satisfy also some symmetry relations. These symmetry relations of the measures $G_r(\10 )$ will imply Theorem B and its generalizations.

It is also possible to prove Conjecture A, as well as its generalization for $\sigma \in G_\bQ$, in this way, 
but the calculations will be much more complicated.

\medskip

\section{The power series associated with Galois action.}

\smallskip
 
Let $V$ be a smooth algebraic variety over a number field $K$ and let $V_{\bar K}:=V\times _K \bar K$. Let $x$ be a point of $V$ with values in an 
algebraic closed field. We denote by
\[
 \pi _1(V_{\bar K},x)
\]
the \'etale fundamental group of $V_{\bar K}$ based at $x$. Assume that $x$ and $y$ are $K$-points or generic points (tangential points) ``defined'' over $K$.
The Galois group $G_K$ acts on $\pi _1(V_{\bar K},x)$ and on the $\pi _1(V_{\bar K},x)$-torsor of \'etale paths on $ V_{\bar K}$ from $x$ to $y$.
Let $\ga$ be an \'etale path on  $ V_{\bar K}$ from $x$ to $y$.
For any $\si \in G_{K}$ we define 
\[
 \ffk _\ga (\si ):= \ga ^{-1}\cdot \si (\ga )\in \pi _1(V_{\bar K},x)\,.
\]
The function $\ffk_\ga :G_K\to \pi _1(V_{\bar K},x)$ is a cocycle and have the following properties:
\begin{enumerate}
 \item[a)] naturality -- if $g:V\to W$ is a smooth algebraic morphism defined over $K$ then
 \begin{equation} \label{eq:*0}
g_*(\ffk _\ga (\si ))=\ffk _{g_*(\ga )}(\si )\,, 
\end{equation}
where $g_*$ is the map induced by $g$ on \'etale fundamental groups and on torsors of paths;
 \item[b)] compatibility with composition of paths -- if $\al$ is a path from $x$ to $y$ and $\be$ from $y$ to $z$, we denote by $\be \cdot \al$ the composed 
 path from $x$ to $z$. Then we have
 \begin{equation} \label{eq:*1}
\ffk _{\be \cdot \al} (\si )=\al ^{-1}\cdot \ffk _{\be}(\si ) \cdot \al  \cdot \ffk _{\al}(\si )\,.
\end{equation}
 \item[c)] Hence we get that 
\begin{equation} \label{eq:*2}
\ffk _{  \al ^{-1}} (\si )=\al  \cdot \ffk _{\al}(\si )^{-1}\cdot \al ^{-1}\,.
\end{equation}
\end{enumerate}

\medskip

We assume that $K\subset \bar \bQ$ and we fix an embedding of $\bar \bQ$ into the field of complex numbers $\bC$.
Hence we have comparison homomorphism
\[
 \pi _1 (V(\bC ),x) \to \pi _1(V_{\bar K},x)\,.
\]
Let us fix a rational prime number $p$ and let us denote 
\[
 \xi _{p^n} :=\exp ({\frac{2\pi i}{p^n}})\,.  
\]
Let us set
\[ 
V_n:=\bP^1_{\bar \bQ}\setminus (\{0,\infty\}\cup \mu _{p^n})\,.                  
\]
Let  $x_n$ (loop around $0$) and $y_{k,n}$ (loop around $\xi ^k_{p^n}$) for $0\leq k<p^n$ be the standard generators of $\pi _1(V_n,\01 )$ (see \cite[page 19]{W3}).
Let 
\[
 \Yc _n:=\{X_n,Y_{0,n},Y_{1,n},\ldots ,Y_{p^n-1,n}\}
\]
and let 
\[
 \bQ _p\{\{\Yc _n\}\}
\]
be a $\bQ _p$-algebra of non-commutative formal power series on elements of $\Yc_n$ and let $I_n$ be its augmentation ideal.
Let
\[
 E_n:\pi _1(V_n,\01 )\to \bQ _p\{\{\Yc _ n\}\}
\]
be a continuous, multiplicative map defined by
\[
 E_n(x_n):=\exp X_n\;\;{\rm and}\;\;E_n(y_{k,n}):=\exp Y_{k,n}\;\;{\rm for}\;\;0\leq k<p^n\,.
\] 
Let $\Mc _n$ be the set of all monomials in non-commuting variables belonging to $\Yc _n$.
If $w\in \Mc _n$ we denote by $\deg w$ the degree of $w$ as a monomial in variables $X_n,Y_{0,n},\ldots ,Y_{p^n-1,n}$ and by $\deg _\Yc w$ the degree of $w$ in variables 
$ Y_{0,n},\ldots ,Y_{p^n-1,n}$.
\medskip

Let $\pi _n$ be the canonical path from $\01$ to ${\frac{1}{p^n}}\10$ ( the interval $[0,1]$). Let $w\in \Mc _n$.
For $\si \in G_{\bQ}$ we define coefficients $\lam _w ^{(n)}$ and $li _w^{(n)}$   
by the equalities
\[
 \Lam _{\pi _n}(\si ) =1+\sum _{w\in \Mc _n} \lam _w ^{(n)}(\si )w\in \bQ _p\{\{\Yc _{n}\}\}
\]
and  
\[
 \log \Lam _{\pi _n}(\si ) = \sum _{w\in \Mc _n} li _w ^{(n)}(\si )w\in \bQ _p\{\{\Yc _{n}\}\}\,.
\]
If $n=0$ we shall usually omit the index $0$ and we shall write 
\[
 \Lam _{\pi  }(\si ) =1+\sum _{w\in \Mc  } \lam _w (\si ) w\in \bQ _p\{\{X,Y \}\}\,.
\]
We shall also omit $\si $ from the formulas to stress that $\lam _w ^{(n)}$ and $li _w^{(n)}$ are functions on $G_{\bQ}$.

\medskip

The study of the coefficients $\lam _w$ and $li _w$ of the power series $\Lam _\pi$ and $\log \Lam _\pi$ is the principal aim of the paper. 
Let $w=YX^{n_1}YX^{n_2}\ldots YX^{n_k}$. The coefficients $\lam _w$ or $li _w$ we view as analogues of the multi-zeta numbers $\zeta (n_1+1,n_2+1,\ldots ,n_k+1)$. For example in 
\cite[Proposition 3.1]{W7} it is shown that 
\[
 \lam _{YX^{2n-1}}=li_{YX^{2n-1}}=-{\frac{B_{2n}}{2\cdot(2n)!}}(\chi ^{2n}-1)\,,
\]
where $B_{2n}$ is the $2n$-th Bernoulli number and $\chi :\Gal (\bQ(\mu _{p^\infty })/\bQ)\to \bZ _p^\times$ is the $p$-cyclotomic character.

\medskip

In the next proposition we present one of our two main tools used in the paper. The proposition below is a special case of Propositions 2.3 and 2.5 in \cite{W9}.

\medskip

\noindent{\bf Proposition 1.1.} Let $\si \in G_{\bQ}$.
\begin{enumerate}
 \item[i)] The family of functions
 \[
 \big\{ G_r^{(n)}(\si ):(\bZ/p^n\bZ)^r\ni (i_1,\ldots ,i_r)\mapsto \lam _{Y_{i_1,n}Y_{i_2,n}\ldots Y_{i_r,n}}^{(n)}(\si )\in \bQ _p\big\} _{n\in \bN}
 \]
forms a measure $G_r (\si )$ on $(\bZ _p)^r$ with values in $\bQ _p$\,.
\item[ii)] Let $w=X^{n_0}YX^{n_1}YX^{n_2}\ldots YX^{n_r}$. Then we have
\[
 \lam _w=\big( \prod _{n=0}^rn_k!\big) ^{-1}\int _{(\bZ_p)^r}(-x_1)^{n_0}(x_1-x_2)^{n_1}\ldots (x_{r-1}-x_r)^{n_{r-1}}x_r^{n_r}dG_r(x_1,\ldots, x_r)\,.
\]

\end{enumerate}

There is an analogous formula for the coefficients  
$li _w$ (see \cite{W9}). In \cite{W9}  the measure  $G_r$  was  denoted by $G_r(\10 )$.

\medskip

\section{The rhombus relation.}

In this section we present our second tool. We start with few definitions. Let $a\in \bP ^1(\bC)$ and let $v$, $ w$ be two tangent vectors at $a$ such that $\mid \mid v\mid \mid=\mid \mid w\mid \mid $.
We denote by $s_a(w,v)$ a path on $\bP ^1 _{\bC}\setminus \{a\}$ from $v$ to $w$ in infinitesimal 
neighbourhood of the point $a$. This path is  an arc in the opposite clockwise sense (see Picture 1).

$$
\,
$$

$$\bigcap _{k=0}^{\infty}\Lc _k^{(n)}
\,
$$

\[
 {\rm Picture\; 1}
\]
Let us define a morphism 
\[
 k_n:V_n\to V_n
\]
by $k_n(\zfk )={\frac{1}{\zfk}}$. Let us set 
\[
 q_n:=k_n (\pi _n)^{-1}\,.
\] 
Let 
\[negatif
 R_{n,1}:V_n\to V_n
\]
be given by $R_{n,1}(\zfk)=\xi _{p^n}\zfk$.

\medskip
  
\noindent
Let us set
\[
 s=s_0( \overset{\to}{0\xi _{p^n}},\01)^{-1},\; t=s_1(\10,\overset{\to}{1\infty})^{-1},\; \eta =R_{n,1}(k_n(s))\; {\rm and}\;
\]
\[
 e=s_{\xi _{p^n}}( \overset{\to}{\xi _{p^n}\infty},\overset{\to}{\xi_{p^n}0})^{-1}\,.
\]
We set also
\[
 c_n=R_{n,1}(\pi _n)^{-1}\;\;{\rm and}\;\; d_n=R_{n,1}(q_n)^{-1}=R_{n,1}(k_n(\pi _n))\,.
\]
Observe that the composition of paths
\begin{equation}\label{eq8}\bigcap _{k=0}^{\infty}\Lc _k^{(n)}
 s\cdot c_n\cdot e\cdot d_n\cdot \eta \cdot q_n\cdot t\cdot \pi _n=1
\end{equation}
in $\pi _1 (V_n,\01 )$ (see \cite[page 167]{W10}, where a similar composition of paths appears).
The picture below has a shape of an octagon or a rhombus, hence a name of a relation we shall deduce.

\medskip

\noindent
Let us set
\[
 \al _1:=t\cdot \pi _n,\;\al _2:=q_n\cdot \al _1,\; \al _3:=\eta \cdot \al _2,\; 
 \al _4:=d_n\cdot \al _3,\; \al _5:=e\cdot \al _4\;{\rm and}\;\al _6:=c_n\cdot \al_5\,.
\]
\medskip

\noindent{\bf Proposition 2.1.} (Octagon relation) On the group $G_{\bQ}$ we have 
\[
 1=\al _6^{-1}\cdot \ffk _s \cdot \al _6 \cdot \al _5^{-1}\cdot \ffk _{c_n} \cdot \al _5 \cdot \al _4^{-1}\cdot \ffk _e \cdot \al _4 \cdot 
 \al _3^{-1}\cdot \ffk _{d_n} \cdot \al _3 \cdot \al _2^{-1}\cdot \ffk _{\eta} \cdot \al _2 \cdot \al _1^{-1}\cdot \ffk _{q_n} \cdot \al _1 \cdot
 \pi _n^{-1}\cdot \ffk _t\cdot \pi _n \cdot \ffk _{\pi _n}\,.
\]

\medskip

\noindent{\bf Proof.} The proposition follows immediately from the formula  \eqref{eq:*1} applied several times to the equality \eqref{eq8}.
\hpb

\noindent{\bf Lemma 2.2.} We have
 \[
 \ffk _{q_n}=q _n^{-1}\cdot (k_n)_*(\ffk _{\pi _n}^{-1})\cdot q_n\,.
\]

\medskip

\noindent{\bf Proof.}
The lemma  follows from the naturality property \eqref{eq:*0} and the formula \eqref{eq:*2}. \hpb

\medskip

\noindent{\bf Lemma 2.3.}  On the subgroup $G_{\bQ (\mu _{p^\infty })}$ of $G_\bQ$ we have 
\begin{enumerate}
 \item[i)] $\ffk _{c_n}=c_n^{-1}\cdot ( R_{n,1})_*(\ffk _{\pi _n}^{-1})\cdot c_n$\,;
 \item[ii)] $\ffk _{d_n}=(R_{n,1})_*((k_n)_*(\ffk _{\pi _n}))$\,;
 \item[iii)] $ \ffk _\eta =(R_{n,1})_*\big( (k_n)_*(\ffk _s)\big)\,.$
 \item[iv)] $\ffk _s=\ffk _e=\ffk _\eta =\ffk _t=1\,.$
\end{enumerate}

\medskip

\noindent{\bf Proof.}  The morphism $R_{n,1}$ commutes with the action of the Galois group $G_{\bQ (\mu _{p^\infty })}$. 
Hence the points i), ii) and iii) of the proposition follow.
Studying the effect of $\ffk _s(\si )$ on the test functions $\zfk^{\frac{1}{p^n}}$ one shows that $\ffk _s(\si )=x_n^{{\frac{1}{p^n}}(1-\chi (\si ))}$
for $\sigma \in G_{\bQ (\mu _{p^n})}$. 
Hence it follows that $\ffk _s(\si )=1$ for $\si \in G_{\bQ (\mu _{p^\infty })}$.
This implies that also $\ffk _e=\ffk _\eta =\ffk _t=1\,.$   \hpb

\medskip

The elements $x_n, y_{0,n},\ldots ,y_{p^n-1,n}$ are free generators of a free pro-$p$ group 

\noindent
$\pi _1 (V_n,\01 )$. Hence the element $ \ffk _{\pi _n}$
is a convergent infinite product of commutators in these generators as the group  $\pi _1 (V_n,\01 )$ is pro-unipotent. We shall write 
$\ffk _{\pi _n}=\ffk _{\pi _n}(x_n, y_{0,n},\ldots ,y_{p^n-1,n})$ to indicate this dependence on  generators. Moreover if 

\noindent
$g:\pi _1 (V_n,\01 )\to G$ is a continuous morphism of groups then
$$g_*(\ffk _{\pi _n}(x_n, y_{0,n},\ldots ,y_{p^n-1,n}))=\ffk _{\pi _n}(g_*(x_n),g_*( y_{0,n}),\ldots ,g_*(y_{p^n-1,n})).$$

\medskip

\noindent{\bf Proposition 2.4.} (Rhombus relation)  We have the following equality on the subgroup $G_{\bQ (\mu _{p^\infty })}$ of $G_\bQ$
\[
 \ffk _{\pi _n}^{-1}\big( \al _6^{-1}\cdot ((R_{n,1})_*(x_n))\cdot \al _6,\al _6^{-1}\cdot ((R_{n,1})_*(y_{0,n}))\cdot \al _6,\ldots ,\al _6^{-1}\cdot ((R_{n,1})_*(y_{p^n-1,n}))\cdot \al _6\big)\cdot 
 \]
 \[
  \ffk _{\pi _n}\big( \al _3^{-1}\cdot ((R_{n,1}\circ k_n)_*(x_n))\cdot \al _3,\al _3^{-1}\cdot ((R_{n,1}\circ k_n)_*(y_{0,n}))\cdot \al _3,\ldots \big)\cdot 
\]
\[
  \ffk _{\pi _n}^{-1}\big( \al _2^{-1}\cdot ((k_n)_*(x_n))\cdot \al _2,\al _2^{-1}\cdot ((k_n)_*(y_{0,n}))\cdot \al _2,\ldots \big)\cdot
  \ffk _{\pi _n}(x_n,y_{0,n},y_{1,n},\ldots )=1\,.
\]

\medskip

\noindent{\bf Proof.} The proposition follows from Proposition 2.1 and Lemmas 2.2 and 2.3.
\hpb

\medskip

Now we shall describe the maps induced by $k_n$ and $R_{n,1}$ on fundamental groups. For our purpose the following result will be sufficient.

\medskip

\noindent{\bf Lemma 2.5.} We have the following equalities and congruences in the group
$\pi _1 (V_n,\01 )$ modulo the commutator subgroup $\big( \pi _1 (V_n,\01 ),\pi _1 (V_n,\01 )\big)$.
\begin{enumerate}
 \item [i)] \[\al _2^{-1}\cdot (k_n)_*(x_n))\cdot \al _2=y_{1,n}^{-1}\cdot y_{2,n}^{-1}\cdot y_{p^n-1,n}^{-1}\cdot x_n^{-1}\cdot y_{0,n}^{-1}\,,\]
 \[
  \al _2^{-1}\cdot (k_n)_*(y_{i,n})\cdot \al _2\equiv y_{p^n-i,n} \;\;{\rm for}\;\; 0<i<p^n\,,
 \]
\[
 \al _2^{-1}\cdot (k_n)_*(y_{0,n})\cdot \al _2= y_{0,n}\,,
\]
 \item [ii)] 
 \[
  s\cdot ( R_{n,1}) _* (x_n)\cdot s^{-1} \equiv x_n\,, \;\; s\cdot  ( R_{n,1})_* (y_{i,n})\cdot s^{-1}\equiv y_{i+1,n}\;\;{\rm for}\;\; 0\leq i<p^n\,.
 \]
\end{enumerate}

\medskip

\noindent{\bf Proof.} The proof we left as an exercise.
\hpb

\medskip

\section{Filtrations of the group Gal$(\bar \bQ /\bQ)$.}

With the action of the Galois group $G_{\bQ}$ on fundamental groups one can associates several filtrations of $G_{\bQ}$. Let us set 
\[
 \Lc _0:=G_{\bQ(\mu _{p^\infty})} 
\]
and 
\[
 \Lc _k:=\{ \si \in \Lc _0 \mid \, \forall n\in \bN \, \forall w\in \Mc _n, 
{\rm deg}w={\rm deg}_{\Yc} w\;{\rm and} \;{\rm deg}w\leq k\,\Rightarrow \lam _w^{(n)}(\si )=0\}
\]
for $k>0$.

\medskip

\noindent{\bf Lemma 3.1.} Let $r$ be a positive integer. Let 
\[
w=Y_{i_1,n} X_n^{n_1} Y_{i_2,n }X_n^{n_2}\ldots X_n^{n_{r-1 }} Y_{i_r,n} X_n^{n_r  } \in \Mc _n
\]
be such that 
$r=\deg _{\Yc}w\leq k$. Then $\lam _w^{(n)}$ is zero on $\Lc _k$.

\medskip

\noindent{\bf Proof.} It follows from the proof of \cite[Proposition 2.6.]{W9} that $ \lam _w^{(n)}=$
\[
 (\prod _{i=0}^ rn_i!)^{-1}\int _{(i_1,\ldots i_r)+(\bZ_p)^r}(-x_1)^{n_0}(x_1-x_2)^{n_1}\ldots (x_{r-1}-x_r)^{n_{r-1}}x_r^{n_r}dG_r(x_1,\ldots, x_r)\,.
\]
The assumption that $\si \in \Lc _k$ implies that the measures $G_r(\si )$ are zero measures for $r\leq k$. Hence the coefficients $\lam _w^{(n)}$ vanish on $\Lc _k$.
\hpb 

\medskip

Let us fix $n$. Now we shall define two filtrations associated with the action of $G_{\bQ}$ on $\pi _1(V_n,\01 )$.
Let us set
\[
 \Lc _0^{(n)}:=G_{\bQ(\mu _{p^\infty })}
\]
and 
\[
 \Lc _k^{(n)}:=\{ \si \in G_{\bQ(\mu _{p^\infty })}\mid \forall w\in \Mc _n,\, \deg _\Yc w\leq k \Rightarrow \lam _w^{(n)}(\si )=0 \}
\]
for $k\geq 1$. The filtration $\{\Lc _k^{(n)}\} _{k\in \bN}$ of $G_\bQ$ is usually called the depth filtration associated with the action of $G_{\bQ}$ on $\pi _1(V_n,\01 )_{pro-p}$.

\medskip

\noindent{\bf Lemma 3.2.} We have
\[
 \Lc _k=\bigcap _{n=0}^{\infty}\Lc _k^{(n)}\,.
\]
\medskip

\noindent{\bf Proof.} It follows from Lemma 3.1 that $\Lc _k\subset \Lc _k^{(n)}$. Hence $\Lc _k\subset \bigcap _{n=0}^{\infty}\Lc _k^{(n)}$. The inclusion of
$\bigcap _{n=0}^{\infty}\Lc _k^{(n)} $ into $\Lc _k$ follows immediately from the definition of the filtrations. \hpb

\medskip

Let us set
\[
 \Gc_0^{(n)}:=G_{\bQ(\mu _{p^\infty })}
\]
and
\[
 \Gc_k^{(n)}:=\{ \si \in \Gc_{k-1}^{(n)}\mid \forall w\in \Mc _n,\, \deg w=k \Rightarrow \lam _w^{(n)}(\si )=0 \}
\]
for $k>1$. The filtration $\{\Gc _k^{(n)}\}_{k\in \bN}$ of $G_\bQ$ is associated with the action of $G_\bQ$ on the quotients of  $\pi _1(V_n,\01 )$
by the terms of the lower central series.

\medskip

\noindent{\bf Lemma 3.3.} For any $n\geq 0$ we have
\[
 \bigcap _{k=0}^{\infty}\Lc _k^{(n)}=\bigcap _{k=0}^{\infty}\Gc _k^{(n)}\,.
\]
\medskip

\noindent{\bf Proof.} Observe that
\[
 \bigcap _{k=0}^{\infty}\Gc _k^{(n)}=\{\si \in G_{\bQ(\mu _{p^\infty })}\mid \forall w\in \Mc _n,\,  \lam _w^{(n)}(\si )=0\}\,.
\]
But this is also $\bigcap _{k=0}^{\infty}\Lc _k^{(n)} $. \hpb

\medskip

\section{Symmetries of the measures $ G_r$. }

Let us denote the summation $\sum _{0\leq i_1,\ldots ,i_r<p^n}$ by $\sum _{r,n}$.
The principal result of this section is the following theorem.

\medskip 

\noindent{\bf Theorem 4.1.} Let $\si \in \Lc _{r-1}$. Then
 \[G_r(x_1,x_2,\ldots ,x_r)(\si )-G_r(-x_1,-x_2,\ldots ,-x_r)(\si )+
 \]
 \[
  G_r(-x_1+1,-x_2+1,\ldots ,-x_r+1)(\si )-G_r(x_1-1,x_2-1,\ldots ,x_r-1)(\si )=0\,.\]

\medskip

\noindent{\bf Proof.} To simplify the notation let us set
\[
 a^n_{i_1,\ldots ,i_r}:=\lam _w^{(n)}\,,
\]
where $w=Y_{i_1,n}Y_{i_2,n}\ldots Y_{i_r,n}$. Let $\si \in \Lc _{r-1}$. Then it follows from Lemma 3.1  that
\[
 \Lam _{\pi _n}(\si )\equiv 1+\sum _{0\leq i_1,\ldots ,i_r<p^n}a^n_{i_1,\ldots ,i_r}(\si )Y_{i_1,n} \ldots Y_{i_r,n}\;\;{\rm modulo}\;\; I^{r+1}_n\,.
\]
It follows from the rhombus relation and from Lemmas 2.3 and 2.5 that
\[
 \big( 1-\sum _n a^n_{i_1,\ldots ,i_r}(\si )Y_{i_1+1,n} \ldots Y_{i_r+1,n}\big)\cdot
 \big( 1+\sum _n a^n_{i_1,\ldots ,i_r}(\si )Y_{-i_1+1,n} \ldots Y_{-i_r+1,n}\big)\cdot
 \]
 \[
 \big( 1-\sum _n a^n_{i_1,\ldots ,i_r}(\si )Y_{-i_1,n} \ldots Y_{-i_r,n}\big)\cdot
 \big( 1+\sum _n a^n_{i_1,\ldots ,i_r}(\si )Y_{i_1,n} \ldots Y_{i_r,n}\big)\equiv 1
\]
modulo $I^{r+1}_n$. Comparing coefficients at monomials $Y_{i_1,n} \ldots Y_{i_r,n}$ we get the identity of measures.  \hpb

\medskip

\section{Euler relations.}

Now we shall formulate and prove our main result. We recall from Introduction and from section 1 that
\[
 \lam _{YX^{n_1}YX^{n_2}\ldots X^{n_r-1}YX^{n_r}}(\si )
\]
are the coefficients of the power series $\Lam _\pi (\si )\in \bQ \{\{X,Y\}\}$, 
where $\pi$ is the canonical path from $\01$ to $\10$ on $\bP^1_{\bar \bQ}\setminus \{0,1,\infty \}$.
\medskip

\noindent{\bf Theorem 5.1.} Let $\si \in \Lc _{r-1}$ and let 
$w=Y  X ^{n_1} Y X ^{n_2}\ldots X^{n_{r-1 }} Y X^{n_r  }$.
If $\sum _{i=1}^rn_i$ is odd then 
\[
 \lam _w (\si )=0\,.
\]

\medskip

\noindent{\bf Proof.}  
Let us set
\[
 d\mu (x_1,\ldots ,x_r)(\si ):=d \big( G_r (x_1,\ldots ,x_r)(\si )-
 G_r (-x_1,\ldots ,-x_r)(\si )+
\]
\[
 G_r (-x_1+1,\ldots ,-x_r+1)(\si )-G_r (x_1-1,\ldots ,x_r-1)(\si )\big)\,.
\]
Let $\si \in \Lc _{r-1}$. Let $n_1,\ldots ,  n_{r-1},q$ be any sequence of length $r$  of non negative integers.
It follows from Theorem 4.1 that 
\begin{equation}\label{eq:**1}
\int _{(\bZ_p)^r}(x_1-x_2)^{n_1}\ldots (x_{r-1}-x_r)^{n_{r-1}}x_r^q \,d \mu (x_1,\ldots ,x_r)(\si ) =0\,.
\end{equation} 
Let us set $m=\sum _{i=1}^ {r-1}n_i$. 
Let us define
\[
 P_q(x_r):=  x_r^q-(-1)^{m+q}x_r^q+(-1)^{m+q}(x_r-1)^q-(x_r+1)^q \,.
\]
After  changes of variables in   the last three integrals in the formula \eqref{eq:**1} we get
\begin{equation}\label{eq:**2}
 \int _{ (\bZ_p)^r}(x_1-x_2)^{n_1}\ldots (x_{r-1}-x_r)^{n_{r-1}}  P_q(x_r) dG_r(x_1,\ldots ,x_r)=0\,.
\end{equation}

\medskip

We shall prove the theorem by induction with respect to $a$   in monomials

\noindent
$Y  X ^{n_1} Y X ^{n_2}\ldots X ^{n_{r-1 }} Y  X ^{a  }$.
Observe that 
\[
P_2(x_r)=   \left\{
\begin{array}{ll}
-4x_r & \text{ if } m \text{ is even}, \\
-2 & \text{ if }  m \text { is odd}\;.
\end{array}
\right.
\]
The identity \eqref{eq:**2}, which holds for any sequence of length $r$ of non negative integers $n_1,\ldots ,n_{r-1},q$ implies that if $m$ is odd
and $q=2$ then
\[
 \int _{ (\bZ_p)^r}(x_1-x_2)^{n_1}\ldots (x_{r-1}-x_r)^{n_{r-1}} \cdot (-2)\, dG_r (x_1,.. ,x_r)=0
\]
on $\Lc _{r-1}$
and if $m$ is even and $q=2$ then
\[
 \int _{ (\bZ_p)^r}(x_1-x_2)^{n_1}\ldots (x_{r-1}-x_r)^{n_{r-1}} \cdot (-4x_r)\, dG_r (x_1,.. ,x_r)=0
\]
on $\Lc _{r-1}$.   Hence the theorem holds for $a=0$ as then $m+0=m$ is odd and for $a=1$ as then $m+1$ is odd.

Let us assume that the theorem is true for $k<a-1$, i.e., if $k<a-1$ and $m+k$ is odd then 
\[
 \int _{(\bZ_p)^r}(x_1-x_2)^{n_1}\ldots (x_{r-1}-x_r)^{n_{r-1}}x_r^k \,d G_r (x_1,\ldots ,x_r)(\si )=0
\]
for $\si \in \ \Lc _{r-1}$.
Observe that 
the polynomial $$P_a(x_r)=\sum _{i=1}^a{a \choose a-i}((-1)^{m-(a-i)}-1)x_r^{a-i}\,.$$
It follows from the equality \eqref{eq:**2} for $q=a$  that
\begin{equation}\label{eq:6}
\sum _{i=1}^a
\int _{(\bZ_p)^r}
(x_1-x_2)^{n_1}\ldots (x_{r-1}-x_r)^{n_{r-1}}\cdot 
\end{equation}
\[
{a \choose a-i}\big( (-1)^{m-(a-i)}-1\big) x_r^{a-i}
dG_r (x_1,\ldots ,x_r)=0\,.
\]
 Notice that the equation \eqref{eq:6} holds for any $m=\sum _{i=1}^{r-1}$ and any $a$.

Let us assume that
\[
 \sum _{i=1}^{r-1}n_i+(a-1)=m+(a-1)
\]
is odd. Then $m+a$ is even. Therefore $(-1)^{m-(a-i)}-1=0$ for $1\leq i \leq a$ and   $i$ even.
For $1\leq i\leq a$ and   $i$ odd we have 
\[
 (-1)^{m-(a-i)}-1=-2\,.
\]
The integrals 
\[
\int _{ (\bZ_p)^r}
(x_1-x_2)^{n_1}\ldots (x_{r-1}-x_r)^{n_{r-1}}\cdot 
{a \choose a-i}(-2)x_r^{a-i}
dG_r (x_1,\ldots ,x_r)
\]
vanish for $i=3,5, \ldots$ and $i\leq a$ by the inductive hypothesis, as then $a-i<a-1$ and $\sum _{i=1}^{r-1} n_i+a-i$ is odd.
Therefore it follows from \eqref{eq:6} that the integral
\[
\int _{ (\bZ_p)^r}negatif
(x_1-x_2)^{n_1}... (x_{r-1}-x_r)^{n_{r-1}}\cdot 
{a \choose a-1}(-2)x_r^{a-1}
dG_r (x_1,.. ,x_r)=0\,.
\] 
Hence we have shown the statement for $a-1$. \hpb 

\medskip
 
\noindent 
{\bf Proof of Theorem B.}
Observe that 
\[
 \Lc _1=\{\si \in G_{\bQ (\mu _{p^\infty })}\mid \forall n \in \bN \, \forall 0\leq i<p^n, \lam ^{(n)}_{Y_{i,n}}(\si )=0\}.
\]
The functions $\lam ^{(n)}_{Y_{i,n}}:G_{\bQ (\mu _{p^\infty })} \to \bZ _p$ are Kummer characters $\kappa (1-\xi ^i_{p^n})$ associated with $1-\xi ^i_{p^n}$ for $n\in \bN$ and $0\leq i<p^n$. 
Hence they vanish on the Galois group $G_{\Kc _1}$ of the field $\Kc _1$. Therefore we have $G_{\Kc _1}\subset \Lc _1$. 
Hence Theorem B follows from Theorem 5.1. \hpb
 
\medskip

\noindent{\bf Corollary 5.2.} Let $r$ be a positive integer. Let $n_1,\ldots ,  n_{r-1},n_r$ be any sequence of length $r$ of  non negative integers.
Let $m=n_1+\ldots + n_{r}$.
Let $\si \in \Lc _{r-1}$ and let $n \geq0$. Then we have
\[
  \int _{(i_1,.. ,i_r)+p^n(\bZ_p)^r }(x_1-x_2)^{n_1}... (x_{r-1}-x_r)^{n_{r-1}}x_r^{n_r}\, dG_r(x_1,.. ,x_r)(\si )+
\]
\[
 (-1)^{m+1}\int _{(-i_1,.. ,-i_r)+p^n(\bZ_p)^r }(x_1-x_2)^{n_1}... (x_{r-1}-x_r)^{n_{r-1}}x_r^{n_r}\, dG_r(x_1,.. ,x_r)(\si )+
\]
\[
 (-1)^{m}\int _{(-i_1+1,.. ,-i_r+1)+p^n(\bZ_p)^r }(x_1-x_2)^{n_1}... (x_{r-1}-x_r)^{n_{r-1}}(x_r-1)^{n_r}\, dG_r(x_1,.. ,x_r)(\si )+
\]
\[
 -\int _{(i_1-1,.. ,i_r-1)+p^n(\bZ_p)^r }(x_1-x_2)^{n_1}...(x_{r-1}-x_r)^{n_{r-1}}(x_r+1)^{n_r}\, dG_r(x_1,.. ,x_r)(\si )=0\,.
\]

\medskip

\noindent{\bf Proof.} We calculate integrals over the set  $(i_1,.. ,i_r)+p^n(\bZ_p)^r $ against the measure $\mu$ from the proof of Theorem 5.1. 
After changes of variables we get the result. \hpb

Below we shall rewrite the formula from Corollary 5.2 in terms of coefficients $\lam _w^{(n)}$.

\medskip

\noindent{\bf Corollary 5.3.} Let $n_1,\ldots ,  n_{r-1},{n_r}$ be any sequence of length $r$  non negative integers. Let $m=n_1+\ldots + n_{r}$.
Let $\si \in \Lc _{r-1}\cap \Gc_{m-1}^{(N)}$ and let $n \geq0$. 
Let $w=
Y_{i_1,n} X_n^{n_1} Y_{i_2,n }X_n^{n_2}\ldots X_n^{n_{r-1 }} Y_{i_r,n} X_n^{{n_r}  }$,

\noindent
$w_1=Y_{-i_1,n} X_n^{n_1} Y_{-i_2,n }X_n^{n_2}\ldots X_n^{n_{r-1 }} Y_{-i_r,n} X_n^{{n_r}  }$,

\noindent
$w_2=Y_{-i_1+1,n} X_n^{n_1} Y_{-i_2+1,n }X_n^{n_2}\ldots X_n^{n_{r-1 }} Y_{-i_r+1,n} X_n^{{n_r}  }$ and

\noindent
$w_3=Y_{i_1-1,n} X_n^{n_1} Y_{i_2-1,n }X_n^{n_2}\ldots X_n^{n_{r-1 }} Y_{i_r-1,n} X_n^{{n_r}  }$.
Then we have
\[
 \lam _w^{(n)}(\si )+(-1)^{m+1}\lam _{w_1}^{(n)}(\si )+(-1)^{m }\lam _{w_2}^{(n)}(\si )-\lam _{w_3}^{(n)}(\si )=0.
\]

\vglue 2cmnegatif

\vglue 1cm

\noindent Universit\'e de Nice-Sophia Antipolis

\noindent D\'epartement de Math\'ematiques

\noindent Laboratoire Jean Alexandre Dieudonn\'e

\noindent U.R.A. au C.N.R.S., N$^{\rm  o}$ 168

\noindent Parc Valroscountere -- B.P. N$^{\rm  o}$ 71

\noindent 06108 Nice Cedex 2, France

\smallskip

\noindent {\it E-mail address} wojtkow@math.unice.fr

\noindent {\it Fax number} 04 93 51 79 74
\medskip

\noindent{\bf 
\medskip

\end{document}